\documentclass{yearbook20}
\usepackage{tabularx}
\usepackage{booktabs}
\usepackage{graphicx}
\def\DDD{\mathcal{D}}
\def\HHH{\mathcal{H}}
\def\III{\mathcal{I}}

\def\KKK{\mathcal{K}}
\def\LLL{\mathcal{L}}
\def\MMM{\mathcal{M}}
\def\NNN{\mathcal{N}}
\def\TTT{\mathcal{T}}

\def\vecb#1{\boldsymbol{#1}}
\def\iffdef{\, {:\!\iff}\,}
\def\imp{\Rightarrow}
\def\nat{\mathbb{N}}
\def\eps{\,\epsilon\,}
\def\epsi{\,\epsilon_i\,}

\def\validb#1#2{\models_{\ll}#1[#2]}
\def\validbn{\models_{\ll}}

\def\pmap{\stackrel{p}{\mapsto}}
\def\smap{\stackrel{s}{\mapsto}}
\def\emb#1#2{emb_{#1}^{#2}}

\def\y#1{\mathsf{#1}}
\def\yx#1{\mathsf{x_{#1}}}
\def\yA{\mathsf{\Phi}}

\begin{document}

\LongTitle{Infinity is Not a Size}
% The name of your contributed article.
\ShortTitle{Infinity is Not a Size}
% A possibly shorter name to appear in the header. This command needs not be used and may be removed if your article title is short enough.

\AuthorA{Matthias Eberl}{Independent Researcher}{Germany}{matthias.eberl@mail.de}
%\AuthorA{}{}{}{}
% The four obligatory arguments of \AuthorA are name, affiliation, state, and e-mail. Affiliation should be ordered from general (e.g., university) to specific (e.g., faculty, or department). If you do not have an affiliation, leave the second argument empty.
%\AuthorAThanks{First author acknowledgement and thanks.}
% Acknowledgement is not obligatory and may be removed.

\Abstract{We present a dynamic model theory that avoids the paradoxes stemming from completed infinities, but does not require any translation of formulae. The main adoption is the replacement of an actual infinite carrier set by a potential infinite one, and by a finitistic interpretation of the universal quantifier.}
% A short summary of your article. The \Abstract command is obligatory.

% A list of a few important keywords in your article. This command is obligatory as well.

\MakeTitlePage
% This command builds up the title page, header and footer, and sends the information about authors at the end of the document. Do not change its position (i.e. it has to appear right after \Keywords) or the document might fail to build.

\section{Introduction}

Our aim is to develop mathematics with the notion of potential infinity instead of actual infinity. Most of the technical details have been worked out and submitted partly. In this paper we give a motivation and general overview of this approach.

\subsection{The Idea in a Nutshell}

The sole and simple starting point is to regard infinity as a potential infinite. The subsequent concepts follow from that in an (almost) natural way. First, the potential infinite is basically a dynamic conception: An infinite set is a process that exhausts its elements, it is not a ``flat'' totality of objects. Its fundamental property is the ``indefinite extensibility of the finite''. The possibility to always extend a totality of objects is more fundamental than (and in contrast to) having a completed totality. Within that process, indefinitely large finite states are substitutes for the actual infinite set.

Second, thinking of sets as processes also requires to take the dependency of different process states into account. For instance, the state of a model may depend on that of the syntax, or the state of a subset depends on the state of the whole set. In this way, the paradoxes of the infinite are avoided, see Section \ref{remdepsec}.

Third, a potential infinite set requires a non-tautological reading of the universal quantifier (Section \ref{reintunivsec}). We cannot simply interpret ``for all'' by referring to a completed infinite set. It is crucial that any reference to an infinite set of elements requires a reference to some specific stage; otherwise one would introduce actual infinity through the back door or already presuppose this interpretation implicitly. In order to formulate this interpretation properly we use a refinement of the Tarskian model theory.

\subsection{Related Work}

In \citep{lavine2009understanding}, in the Section ``The Finite Mathematics of Indefinitely Large Size'' the author introduces the concept of an \emph{indefinitely large} size, which is nevertheless finite. His ideas are based on a work of \citet{Mycielski1986} about locally finite theories. The central idea of Mycielski's work is that the range of the bound variables inside a single formula increases. For instance, in a formula $\forall \yx{0} \exists \yx{1} \yA$, the variable $\yx{0}$ refers to objects in some carrier set $\MMM_{i_0}$ and a further quantification, e.g.~by $\exists \yx{1}$, refers to a set $\MMM_{i_1}$, being indefinitely large relative to $\MMM_{i_0}$. Lavine uses Mycielski's results in order to argue that one is able to extrapolate results from the finite to the infinite\footnote{A critical position on this extrapolation can be found in \citep{sereno1998infinity}. For a critical view on the concept of an ``indefinitely large finite'' see \citep{bremer2007varieties}.}. Both, Mycielski and Lavine, do not investigate the notion of a potential infinite in this context.

Mycielski constructs a finite model for a finite set of formulas, implicitly leading to a dependency of the model on the syntax\footnote{Note that we regard this dependency of the model on the syntax as one of the basic property of the potential infinite.}. In order to interpret formulas of a first-order theory one has to translate them: Each bound variable is restricted by a predicate $\Omega_i$, indicating some upper bound. Additionally one adds bookkeeping axioms in order to formulate the exchangeability of these $\Omega_i$. Thereby Mycielski uses the common Tarskian model theory (for his purpose, a finite term model suffices).

To cope with indefinite extensibility we give a new interpretation of the universal quantifier, presented in Section \ref{logicsec}, related to Mycielski's idea of the increasing ranges of quantifiers. In this way, no translation of formulas is necessary and no further axioms are required. This new model theory also allows the review of meta-mathematical properties, see Section \ref{metasec}, and proposes a way to extend the concepts to higher-order logic, see Section \ref{highersec}.

\citet{lavine2009understanding} calls the translated, finitistic version of the axiom of infinity the ``axiom of zillions''. For him, the axiom of infinity is the extrapolation of the axiom of zillions from finite to infinitary set theory. In our approach ZFC is just another first-order theory. The axiom of infinity claims the existence of an infinite object, which is a representative of a potential infinite set (see Section \ref{highersec}). This set is potential infinite and not actual infinite due to the non-trivial, finitistic interpretation of the universal quantifier. There are no two versions of the axioms (a finitary and an infinitary one) in our approach.

\citet{shapiro2006all} considered the potential infinite as an indefinite extensible concept. Therein they state that ``If a `collection' is not a set, then it is nothing, has no size at all, and so can't be {`too big'}'' and ask: ``The question, simply, is whether it is ever appropriate or intelligible to speak of all of the items that fall under a given indefinitely extensible concept''. The authors come to the conclusion that there is no satisfying solution how to read such a quantification and refer to \emph{reflection principles} as a possible answer. The interpretation that we present has an implicit reflection principle.

To sum up: If one regards the infinite as a potential infinite, one does not dispel the infinite as Mycielski does. And there is no need to use knowledge about finitely many objects in order to get knowledge about infinitely many objects, as Lavine argues. The infinite appears in a different form, it is no longer a size. An infinite set is a process and infinity creates a dependency of several process states, not of absolute cardinalities. Its new basic property is indefinite extensibility and not ultimate completion.

\subsection{What it is Not}

Constructive logic\footnote{Note that intuitionism relies on the notion of a potential infinite, but does not explain it. It refers in a naive way to infinite totalities without using some finite state.} and recursion theory essentially uses the idea of a decision procedure, which plays no role in our consideration. In particular, we do not claim that there is an effective procedure to determine an indefinitely large finite state. Also, to claim the existence of an object with a specific property does not presuppose a procedure that terminates within known bounds. 

We do not need any change of the axiom system as originally done by Jan Mycielski, or in Marcin Mostowskis work (see e.g.~\citealp{mostowski2005fm}). Mostowski uses potential infinite sets $(\nat_i)_{i \in \nat}$ with $\nat_i := \{0, \dots, i-1\}$ as a basis of his model theory, but his approach is less dynamic and the axioms of Peano arithmetic must be adopted in a way that there exists a greatest number (see~\citealp{mostowski2003representing}). We also do not use paraconsistent logics to cope with the challenges of a finitistic view (cf.~\citealp{van1994strict}, \citealp{bremer2007varieties}, \citealp{priest2013indefinite}). Ideas based on ``feasibility'' are not relevant for our approach, an idea formulated e.g.~by \citet{parikh1971existence}. And we do not introduce modalities: Recently \citet{linnebo2019actual} suggested to formalize the potential infinite using a modal reading of this notion. We also avoid vagueness or a notion of a ``grey zone''.

We do not want to introduce or use a philosophical position. Our terminology and the technical treatment does not presuppose any view of mathematics such as Platonism, formalism or intuitionism. If we speak of objects that are created at some stages this should not be understood literally. We might also say that they have been conceived or revealed. Similarly, if we speak of time, later states etc., this again should not be understood literally, but figuratively.

\section{The Potential Infinite}

The locution ``potential infinite'' is a technical term, there is no infinity (as opposed to a finite) involved in this concept. Seen in that way, the potential infinite is a form of finitism, since at each stage there are only finitely many objects. If we refer to a potential infinite set, we necessarily have to refer to some state with finitely many elements. Often finitism is regarded as a view that postulates a \emph{fixed} bound. With the notion of a potential infinite, the common finite vs.~infinite distinction becomes a fixed vs.~variable difference.

\subsection{Ontological and Epistemological Finitism}

The form of finitism that results from the potential infinite is also an \emph{ontological} one, not an \emph{epistemological} one. Another direction of finitism is the restriction of inferences to ``finitary'' ones leading to a ``finitary reasoning''. This is often related to Hilbert's program, which is concerned with the justification of classical mathematics by finitary reasoning. \citet{tait1981finitism} claimed that finitist arithmetic coincides with primitive recursive arithmetic PRA. But his concept of finitism is based on a different idea than those of an infinite totality. As he states in \citep{tait2002remarks}: ``My argument is that one can understand the idea of an arbitrary object of a given finitist type independently of that of the totality of objects of that type''.

This (ontological) finitism is a justifiable position. For instance, \citet{fletcher1989truth} distinguishes between abstraction and idealization on the one hand, being basically a simplification, and extrapolation on the other hand, being an extension to a larger world. He argues that abstraction and idealization are necessary steps to obtain mathematical objects from physical objects. These procedures justify for instance the step from feasibilism to finitism. But the step to an actual infinity requires an extrapolation. Following this argument, a finitistic position is a natural result of the usual abstraction and idealization process from our experiences. The concept of an actual infinite however goes beyond these abstractions and idealizations.

\subsection{Infinite Entities}

There are different infinite entities such as sets, spaces, lines, decimal representations, numbers etc., some numbers are also said to be infinitesimal small. Basically we distinguish between collections of objects --- including structured collections, relations and functions --- and single objects. 

All collections of objects have some underlying set, e.g.~a space has an underlying set of points or a function $f$ an underlying set of assignments $a \mapsto f(a)$, and it suffices to explain when this set is infinite. If we call a single object infinite, e.g.~a number is infinitely large or infinitesimal small, then it has an \emph{infinite set of defining properties} such as $\omega > n$ or $\epsilon < \frac{1}{n}$ for all (infinitely many) natural numbers $n$. Often an infinite object has infinitely many better and better approximations. So first and primarily we have to answer the following question: \emph{When is a set of entities called infinite?}

\subsection{Paradoxes of the Actual Infinite}

The existence of actual infinite sets are based on the ``domain principle'' \citep{hallett1984cantorian}, that every potential infinity presupposes an actual infinity. Indefinite extensibility, as we understand it, is contrary to this principle. Nevertheless, statements still have determined truth values.

Actual infinities have several counter-intuitive properties. This starts with simple examples, e.g.~there are as many natural numbers as even numbers. More complex examples are the Banach-Tarski paradox, as a consequence of the idea that a continuum is an actual infinite set of points. The deficiency of this view is not the fact that actual infinite sets have unfamiliar properties, but the fact that these ``properties'' stem from relations and dependencies between infinite sets which have been removed --- we show this in an exemplary way in Section \ref{remdepsec}. Simply taking these dependencies into account prevents these paradoxes, which arise as self-made problems that have nothing to do with the mathematical content.

And even more, to introduce actual infinities does not eliminate the phenomenon of indefinite extensibility. After establishing infinite sizes in form of ordinal and cardinal numbers, the question arises naturally, what is the size of the totality of these infinite numbers. It is well known that this again leads to contradictions or paradoxes and a satisfying solution is not available\footnote{Some attempts are done with reflection principles (e.g.~that from Lévy Montague) or with hierarchies of Grothendieck universes. Often a notion of small versus large is introduced, being basically a variant of the original distinction of finite versus infinite.}.

\subsection{Removing a Dependency}
\label{remdepsec}

The cause for the paradoxes of infinity is that actual infinite sets do not allow to consider dependencies between the way its elements are exhausted. Potential infinite sets have exactly this additional structure. Each (potential) infinite set has its own way of exhausting its elements, so these processes can be related to each other. If we switch to an absolute completion of this process, not only a temporary stage, this removes the dynamic and dependency. This is best seen at one of the simplest paradoxes: An infinite set has the same size as a proper subset, whereby ``same size'' should mean a one-to-one correspondence of elements. For instance, the set $\nat$ of natural numbers has the same size as the set $2\nat$ of even numbers, given by the bijection $n \mapsto 2n$.

If we regard both sets as processes, there are different ways to relate them. The set $\nat$ cannot be given as a whole, but solely by some state $\nat_i = \{0, 1, \dots, i-1\}$ and similarly $2\nat$, say by the state $2\nat_j = \{0, 2, \dots, 2j-2\}$. The set $2\nat$ is a subset of $\nat$ only if $i \geq 2j - 1$ and there is a one-to-one correspondence only if $i = j$. Both requirements are not meet at the same time if $i > 1$. But if we imagine that $i$ and $j$ arrive at the infinite, given by a state $\omega$, the simultaneous satisfaction of both requirements becomes possible (intuitively since $\omega$ and $2\omega -1$ are equal as cardinal numbers). So the paradox is removed if we do not allow this step to a completed infinite state, which moreover contains an incontinuity: $i \not= 2i - 1$ for $i > 1$ becomes $\omega = 2\omega - 1$ in the limit.

\subsection{Indefinitely Extensible and Indefinitely Large}
\label{indefexsec}

The basic property of a (potential) infinite set $\MMM$ is its \emph{indefinite extensibility}. But quantification requires reference to some state and the naive interpretation of ``for all $\dots$'' cannot be used --- it is meaningful only if there is a completed set of objects. Hence one needs a fixed finite set that is a replacement of the idealized totality of all possible elements. This is an indefinitely large stage $\MMM_i$ within $\MMM$, relative to a context of other states. An indefinitely large finite set is thus a context dependent version of the actual infinite set. Conversely, if the context can be ignored without harm, one may treat the indefinitely large set as an absolute infinite totality.

The notion of indefinite extensibility that we use includes Dummett's understanding\footnote{Note that whereas Dummett concludes that statements quantifying over an indefinitely extensible concept do not follow the laws of classical logic, our model-theoretic approach does not require this.}. In \citep{dummett1994mathematics} he defines: ``An indefinitely extensible concept is one such that, if we can form a definite conception of a totality all of whose members fall under that concept, we can, by reference to that totality, characterize a larger totality of all whose members fall under it.'' The ordinal numbers and sets are a typical example: If we refer to ``all sets'', this creates or reveals a new set and thus the totality of all sets has changed. 

But already the natural numbers form such an indefinitely extensible concept. If we refer to the number of all numbers, then this reference creates a new number. First there is no number, hence the number of numbers is $0$. So we created a first number, namely 0, and the number of numbers is $1$. Henceforth there are the two numbers $0$ and $1$, creating number $2$ and so on.

The notion of an indefinitely large finite could be seen as a \emph{relative infinite}. If $\III$ denotes the set of states or indices, then a relative infinite is a relation $C \ll i$ (or $i \gg C$) between an index $i \in \III$ and a context $C := (i_0, \dots, i_{n-1})$, with $i_0, \dots, i_{n-1} \in \III$, stating that $i$ is \emph{indefinitely large} or, using a more technical notion, \emph{sufficiently large} relative to $C$. 

We can only investigate finitely many objects in a way that we explicitly refer to them. Say these are currently $a_0, \dots, a_{n-1}$. Most often these objects are not fixed but variable ones, taken from some infinite sets. Assume that $a_0, \dots, a_{n-1}$ are (variable) natural numbers, then saying that $a_k$ is a natural number means $a_k \in \nat_{i_k}$ for some state $i_k \in \nat$. So the currently investigated objects, here $a_0, \dots, a_{n-1}$, are always within a context $C = (i_0, \dots, i_{n-1})$.

By seeing infinity as an indefinitely large finite, the infinite is not outside of an indefinitely extensible set, it is a part of it. It is only outside the region that we can reach from the current stage with our current means. The indefinitely large finite sets $\MMM_i$ with $i \gg C$ behaves exactly as actual infinite sets in the current context of investigation. But they are not completed in an absolute way, i.e., if we change the context, an extension could be necessary.

\begin{figure}[ht]
\centering
\includegraphics{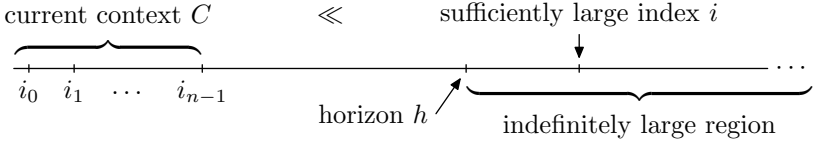}
\caption{The structure of indefinitely extensible sets.}
\label{fig1}
\end{figure}

The notion of an indefinitely large finite is relative in three ways. First, it is not a single state $i \in \III$, but a region, e.g.~$\{i \in \III \mid i \geq h\}$, the \emph{indefinitely large region}. If there is a least element in this region, we call it \emph{horizon}. Secondly, the region depends on a context $C = (i_0, \dots, i_{n-1})$, it is thus a relation $C \ll i$. Figure \ref{fig1}.~illustrates this situation. And thirdly, it is not a single relation $\ll$ but several ones. Their basic properties are that $C \ll i \leq i'$ implies $C \ll i'$ and additionally that $(i_0, \dots, i_{k-1}) \ll i_k$ holds for all $k < n$. The latter expresses a dependency of the size of set $\MMM_{i_k}$ on the sizes of the sets $\MMM_{i_0}, \dots, \MMM_{i_{k-1}}$.

If it is necessary to include the indefinitely large set $\MMM_i$ into the current context $C$ as a further set $\MMM_{i_n}$, then the current context becomes $C' = (i_0, \dots, i_{n-1}, i_n)$. We may then again choose an indefinitely large index $i' \gg C'$.

\section{Mathematical Formalization}

We are now ready to present the mathematical concepts. We start with the basic concept of a system and apply it to first-order logic (FOL). In a second step we treat infinite objects when we switch to higher-order logic (HOL).

\subsection{Systems: Formalizing the Dynamic Aspect}

A dynamic object $a$ is mathematically a function from a index set $\III$ to possible states $a_i$ of $a$, that is, $a$ is the family $(a_i)_{i \in \III}$. Similarly, a dynamic set $\MMM$ is a family of sets varying over $\III$, i.e., $\MMM$ is $(\MMM_i)_{i \in \III}$. The finitistic approach requires that $\MMM_i$ is finite, but the technical treatment makes no use of this property --- we could work with other small versus large distinctions as well (these are often called ``small'' and ``large'', but we try to avoid these notions due to their connotation of size). But from our perspective we gain nothing by doing so. More relevant is the fact that we do not need any notion of finiteness of the sets $\MMM_i$. 

The set $\III$ is equipped with a preorder $\leq$ indicating that if $i' \geq i$, then $i'$ refers to a ``later'' state. We also require that $\III$ is directed with the idea that the different states converge towards ``one limit'' without reaching it\footnote{This ``not reaching the limit'' applies for infinite index sets. It is obvious that any fixed finite set can be subsumed under this concept by considering a singleton index set $\III = \{\ast\}$, or another finite index set with greatest element, which then represents the limit state.}.

The objects in the different states $\MMM_{i'}$ and $\MMM_i$ for $i' \geq i$ are connected by a relation $a_i \smap a_{i'}$ if $a_i \in \MMM_i$ and $a_{i'} \in \MMM_{i'}$. The intention is that $a_{i'}$ is a \emph{successor} of $a_i$. The simplest case for FOL is that $\smap$ is an embedding $\emb{i}{i'} : \MMM_i \to \MMM_{i'}$. In that situation we identify $a_i$ with $\emb{i}{i'}(a_i)$, so $\MMM_i \subseteq \MMM_{i'}$. For instance, $\nat_i \subseteq \nat_{i'}$ and $\nat_i \ni n \smap n \in \nat_{i'}$ for $n < i \leq i'$, which is a special case of a \emph{direct system}. In the sequel we simply speak of a \emph{system} if we refer to $(\MMM_\III, \smap)$.

\subsection{Logic}
\label{logicsec}

The set $\LLL$ of expressions, such as terms and formulas, are typically a free structure over a vocabulary. What is important here is that if $\LLL$ is infinite, it is consequently seen as indefinitely extensible. In all usual languages $\LLL$ one can exhaust its expressions by some measure of complexity, whereby the set $\nat$ as index set suffices\footnote{Consider for instance a first-order language with free variables $\yx{0}, \yx{1}, \dots$, a finite set of relation symbols, the connectives $\neg$, $\land$, $\lor$, $\to$ and the quantifiers $\forall$, $\exists$. Then $\LLL_k$ may consist of the finitely many formulas $\yA(\yx{0}, \dots, \yx{n-1})$ (i.e., free variables within $\yx{0}, \dots, \yx{n-1}$) with $n \leq k$ and less than $k$ connectives and quantifiers.}. So $\LLL = (\LLL_k)_{k \in \nat}$ with finite sets $\LLL_k$ and $\LLL_k \subseteq \LLL_{k'}$ for $k \leq k'$. Each model has an underlying carrier set $\MMM$, probably more. $\MMM$ is formalized as a system $\MMM_\III := (\MMM_i)_{i \in \III}$, for instance $\nat$ is the system\footnote{We sometimes take the index set $\nat^+ := \{1,2,\dots\}$ instead of $\nat$ to avoid empty sets of objects. This is relevant only in HOL, but is an insignificant change.} $(\nat_i)_{i \in \nat}$.

\subsubsection{A Reinterpretation of the Universal Quantifier?}
\label{reintunivsec}

The reinterpretation of the universal quantifier, or better, its finitistic reading, is crucial here. \citet{van1999largest} argues that a reinterpretation of the universal quantifier is indispensable and unassailable for a finitistic point of view:

``In the first place, a classical mathematician or logician will surely remark that, however clever this procedure might be, it still implies a reinterpretation of the universal quantifier. `For all x, $\dots$' does not have its classical meaning, for, in all cases, we are supposed to read `For all $x \leq K, \dots$'. $[\dots]$ if one asks for a standard interpretation of the universal quantifier, then one presupposes the possibility of an infinite domain, hence one can never have such an interpretation in a finite domain.''

A new interpretation makes it impossible to formulate the concept of a potential infinite in the object language. So whether an unbounded set is seen as actual or potential infinite cannot be stated in a formal system, e.g.~as some axiom --- for this reason \citet{niebergall2014assumptions} found no convincing formulation that a theory assumes merely the potential infinite, and not an actual infinite. Only in a formalization of the background theory it is possible to distinguish potential and actual infinity as different interpretations of the universal quantifier.

\subsubsection{Classical First-Order Logic}
\label{clfolsec}

In a FOL the only subject that is considered to be infinite is the domain of discourse and the functions and relations on it, not the objects themselves. Our starting point for classical FOL is the usual Tarskian model theory. The variable assignment $\vecb{a} = (a_0, \dots, a_{n-1})$, instantiating the free variables $\yx{0}, \dots, \yx{n-1}$ in a formula $\yA$, are taken from $\MMM_C := \MMM_{i_0} \times \dots \times \MMM_{i_{n-1}}$ (see Section \ref{indefexsec}). Relation $R$, the interpretation of an $n$-ary relation symbol $\y{R}$ in the signature $\Sigma$, is a family $(R_C)_{C \in \III^n}$ with $R_C \subseteq \MMM_C$. It satisfies $R_C \subseteq R_{C'}$ for $C \leq C'$ with a pointwise order\footnote{It is also possible to add function symbols and functions, too.}.

The notion of validity becomes $\MMM_\III \validb{\yA}{\vecb{a} : C}$ with $\ll$ as an additional parameter. The interpretation of the logical connectives is as usual. The universal quantifier is interpreted by
\begin{equation}
\label{alleq}
\MMM_\III \validb{\forall \yx{n} \yA}{\vecb{a} : C} \ \ \iffdef \ \MMM_\III \validb{\yA}{\vecb{a} b : C i} \text{ for all } b \in \MMM_i
\end{equation}
with an \emph{indefinitely large} set $\MMM_i$, i.e., $i \gg C$. Therein, $\vecb{a} b$ is the extension of the variable assignment $\vecb{a}$ by $b$, and $C i$ the extension of the context $C$ by index $i$. The existential quantifier $\exists \yx{}$ can be reduced to $\neg \forall \yx{} \neg$.

The interpretation thus uses a semantic reflection principle: An assertion about all infinitely many elements in $\MMM_\III$ is true if and only if it is true about a part $\MMM_i$ that reflects the whole infinite set in the current context $C$. In order to show that this interpretation is sound (relative to a usual inference system for classical FOL) and that the validity is independent of the chosen index $i \gg C$, one has to restrict the set of relations $\ll$ to those, which contain all witnesses of valid existential assertions. The construction is similar to the proof of the L\"owenheim-Skolem theorem and Mycielski's original work (see \citealp{Mycielski1986}). The interpretation $\validbn$ for these relations $\ll$ is also \emph{complete}. Moreover it is \emph{locally finite} in the sense that a finite set of expressions requires a finite model only.

\subsubsection{Intuitionistic First-Order Logic}

This interpretation with reflection principle is also applicable to intuitionistic logic with Kripke models, or alternatively Beth models. Kripke models are based on a frame $(\KKK,\leq)$, being a preorder of \emph{epistemic states} or \emph{nodes}. Deviating from the index set $\III$ these frames are typically not directed --- different branches generate alternatives, which may not ``meet later''.

There is however a more important difference between the index sets $\KKK$ and $\III$: The set $\KKK$ represents \emph{epistemic states} or \emph{states of information}, whereas the index set $\III$ is related to the \emph{ontological} side of the semantics. So we may ask whether an element $a$ exists in $\MMM_i$, but if it exists there, it has already all of its properties (nothing new can be added to $a$ in $\MMM_{i'}$ for $i' \geq i$). In comparison, if $a$ exists at $k$, it has the properties known at $k$, but we may discover at $k' > k$ additional properties, even if the same objects exist at $k$ and $k'$. Formally, a relation $R$ in a Kripke structure is a family of sets $R^k_C$ with $k$ a node and $C$ a context. The different requirements become $R^k_C(\vecb{a}) \iff R^{k}_{C'}(\vecb{a})$ for $\vecb{a} \in \MMM_C \cap \MMM_{C'}$ on the one hand and the weaker requirement $R^k_C(\vecb{a})\, \imp\, R^{k'}_C(\vecb{a})$ for $k \leq k'$ and $\vecb{a} \in \MMM_C$ on the other hand.

Similar as for classical logic, one can define a forcing relation with a reflection principle and show soundness and completeness as well as locally finiteness.

\subsubsection{Higher-Order Logic}
\label{highersec}

HOL extends FOL in that it handles infinite objects. A finitistic point of view can deal with infinite objects, such as $\omega$ (the object representing the set of natural numbers in set theory), in two ways: 
\begin{enumerate}
\item The object, say $\omega$, may be regarded as a single abstract object. FOL can deal with these kind of objects only.

\item Or infinite objects are indefinitely extensible objects, analogously as sets are indefinitely extensible. Then there is no single object $\omega$, only approximations of it, e.g.~finite sets $\omega_i$. HOL formalizes this concept as well.
\end{enumerate}

In set theory, as in any other first-order theory, the universe is a system $(\MMM_i)_{i \in \III}$ with finite sets $\MMM_i$. The abstract set-objects $a \in \MMM_i$ become concrete in the background model due to the indefinitely increasing membership relation $\eps$, which is a family with states $\epsi \subseteq \MMM_i \times \MMM_i$. An object such as $\omega$ \emph{represents} an increasing family $\omega_i := \{a \in \MMM_i \mid a \epsi \omega\}$ for all $i \in \III$ with $\omega \in \MMM_i$. Thereby $\omega_i$ is a ``real'' set, i.e., a (finite) set in the background model, that increases if we enlarge $i \in \III$. There is no explicit concept in FOL that sees all of these $\omega_i$ as approximations of $\omega$. The relation $\smap$ in a first-order model is always a function, that is, an object is not differentiated if the model increases (e.g.~in HOL an approximation $1.41$ of a real number differentiates to $1.410, 1.411, \dots 1.419$ if the precision increases from 2 to 3 digits). To put it another way, if a (first-order) object occurs at some stage, it is already ``complete'' and is not an approximation that needs further differentiation.

HOL explicitly introduces approximations which are placed in relation by $\smap$. Infinite objects are only handled by their approximations, e.g.~$\omega$ is the family of finite von Neumann ordinals $n$ and $n \smap n'$ if $n \eps n'$. Then $\omega$ exists only as this family, not as an abstract limit object --- but in HOL there might be both, abstract objects and families of approximations. Moreover, each sufficiently large $n$ plays the role of $\omega$, whereby its size depends on the current context. 

Similarly, a real number, for instance $\sqrt{2}$, can be given as an indefinitely extensible Cauchy sequence of rational numbers (more precisely an equivalence class thereof). Then a sufficiently good approximation by a rational number is a substitute for the real number. Additionally the real number may be given as abstract object of some base type, which is an element of a complete ordered field (defined axiomatically). But $\sqrt{2}$ is \emph{not} the infinite Cauchy sequence, being an element of this field of real numbers.

An elegant way to formulate classical HOL is simple type theory (STT), see e.g.~\citep{farmer2008seven}. A prerequisite for a model is an indefinitely extensible system with a limit construction in the style of a direct and inverse limit. Adding a rule for the universal quantifier (not only a constant) allows an interpretation of the universal quantifier as defined in (\ref{alleq}). 

In STT all infinite objects are formulated as (higher-order) functions. In order to see how STT handles them by approximations, consider functions on natural numbers. Instead of first building an actual infinite set $\nat$ and then defining the function space $[\nat \to \nat]$ as consisting of ``all'' functions, $[\nat \to \nat]$ is approximated by finite function spaces $\NNN_{i \to j} := [\nat_i \to \nat_j]$, with $i, j \in \nat^+$ --- we suggestively write $i \to j$ for such a pair of indices. These function spaces form a system $(\NNN_{i \to j})_{i,j > 0}$ with $f \smap f'$, whereby a function $f : \nat_i \to \nat_j$ is related to $f' : \nat_{i'} \to \nat_{j'}$ for $i \leq i'$ and $j \leq j'$ if $f$ is the restriction of $f'$ to $\nat_i$. A function is itself an indefinitely extensible system $(f_{i \to j})_{(i,j) \in \HHH}$ with $f_{i \to j} : \nat_i \to \nat_j$ and a suitable index set\footnote{For simple functions from a first-order domain to another, e.g.~$f : \nat \to \nat$, it is easy to find these index sets (e.g.~$\HHH = \{i \to j \mid j > \max_{n < i}(f(n))\}$). For functionals of type $2$ or higher, the definition of what is a suitable index set is indeed a challenge.} $\HHH \subseteq \nat^+ \times \nat^+$.

All other constructions to handle infinite objects in a finitistic way use, to our knowledge, completed infinities. A wide spread approach to handle infinite objects very generally is \emph{domain theory}. A domain $\DDD$ is \emph{directed complete}, which is a property that requires actual infinite sets in order to be non-trivial, since each finite directed set has automatically a greatest element. There is another, less well known approach, the hyperfinite type structure, c.f.~\citep{normann1999hyperfinite}. It is based on a Fr\'{e}chet product. In a Fr\'{e}chet product two elements $(a_i)_{i \in \nat}$ and $(b_i)_{i \in \nat}$ are identified if they differ only w.r.t.~finitely many indices $i$. This is obviously trivial if $\nat$ is regarded as potential infinite, since all elements are identified then\footnote{Note that notions such as directed completeness or a Fr\'{e}chet product are meaningful inside a theory if interpreted in indefinitely extensible models. They get a new reading in which the limit is indefinitely large, but not actual infinite. However, we cannot use this understanding to establish the model theory itself, if we do not want to already presuppose this new interpretation.}.

\subsection{Some Meta-Mathematical Considerations}
\label{metasec}

The understanding of the language as increasing is important for meta-mathematical properties. With this understanding some notions collapse, for instance \emph{consistent} is the same as \emph{finitely consistent}. As a consequence, all models are saturated. Interestingly, compactness did not become trivial, but uniformity: If we find a model at each stage of $(\TTT_k)_{k \in \nat}$ (an increasing set of formulas with $\TTT_k \subseteq \LLL_k$), then there is a common model for all stages. This follows from the construction of the Henkin-model in a first-order completeness proof.

We will shortly discuss the presence of non-standard elements in FOL. What is new in our approach is that relation $\smap$ allows a distinction between standard and non-standard models. Systems with partial surjections $\pmap$ (here $\pmap$ is the inverse of $\smap$, indicating a \emph{predecessor} relation) are called \emph{standard}, all others \emph{non-standard}. That is, in a standard model objects are not identified at later states whereas non-standard models allow an \emph{identification} of objects when the elements increase, i.e., $a_0 \smap a$ and $a_1 \smap a$ for $a_0 \not= a_1$ both in $\MMM_i$. HOL additionally allows a \emph{differentiation} in contrast to FOL, i.e., $a \smap a_0$ and $a \smap a_1$ for $a_0 \not= a_1$ both in $\MMM_i$.

Let us consider Peano arithmetic (PA). The non-standard elements arising in Henkin's completeness construction (see e.g.~\citealp{enderton2001mathematical}) are not infinitely large numbers. They are natural numbers that cannot be seen as a number, or for which we do not know, at the current stage of the model construction, which number it is. These non-standard elements are usually considered to be \emph{beyond} all natural numbers $n \in \nat$. A reading of infinite as potential infinite only shows that the number is larger than all finitely many numbers in the current context, which is easily satisfied by a (standard) natural number due to the indefinite extensibility.

Moreover, PA is a theory with equality, so $=$ is interpreted by the identity. In order to satisfy this condition in the common Henkin-construction, a switch to equivalence classes is done \emph{after} all elements have been introduced to the model. The two processes --- adding elements to the model on the one side and identifying them if an equality can be proven on the other --- are seen as independently increasing. One is able to complete the first task (adding elements) before starting the second one (identifying them). 

With a dynamic reading these processes must be done simultaneously in a direct limit construction with non-injective embeddings. So at each stage of the construction there are typically elements (closed terms) that are known to be provably equal to some $\y{n}$ (i.e., the closed term $\y{S} \dots \y{S}\y{0}$ with $n$ successor symbols $\y{S}$). At each stage there will always be some of these unknown elements due to open terms and also due to G\"odel's first incompleteness theorem. In contrast, the usual iterative construction $(\nat_i)_{i \in \nat}$ of a standard model does not introduce non-standard numbers nor identify elements in further steps. Remind that completeness holds with respect to non-standard models (incl.~the standard model), whereas categoricity holds w.r.t.~the standard model only.

\subsection{Application to the Background Model}

The here presented model theory is applicable to its own background model if we apply the idea of indefinitely extensible sets to its meta theory and model. That is, the implicit background model of model theory that we applied in Section \ref{logicsec} makes use of extensible totalities, too. 

In particular, the index set $\III$ is seen as indefinitely extensible (note that $\III$ does not contain objects of the investigated model $\MMM_\III$, but indices are part of the background model). This does not lead to an infinite regress. One might think that since we replaced $\nat$ by $\nat_\nat$, then the index set must be replaced in the same way leading to $\nat_{(\nat_\nat)}$, and so on. But we do not have to perpetually replace $\nat$ by $\nat_\nat$ since there never was a completed set $\nat$ that had to be replaced. From the very beginning there were only indefinitely extensible sets with states $\nat_j$, and this model theory made this explicit. 

In a first reading, the reader may regard an index set $\nat$ as actual infinite. If she accepts the new interpretation, the reader can go through the paper, but now with the new understanding of set $\nat$. Then, if we mention set $\nat_\nat := (\nat_i)_{i \in \nat}$, the index set $\nat$ refers to some stage $\nat_j$ of the background model, i.e., $\nat_\nat$ is $(\nat_i)_{i < j}$. The requirement w.r.t.~the index $j$ is that it must be sufficiently large to describe all investigated models $\MMM_\III$. Therein the locution ``all'' has to be read as ``all indefinitely many'' as defined formally in Formula (\ref{alleq}).

\section{Conclusion}

We presented a natural way how to use a potential infinite in mathematics instead of completed infinities, without any restrictions on the inference system. Infinite sets are seen as indefinitely extensible, realized by an interpretation of the universal quantifier that uses a finitistic reflection principle.

\bibliographystyle{apacite}
\bibliography{InfinityNotSize}

\begin{thebibliography}{}

\bibitem [\protect \citeauthoryear {%
Bremer%
}{%
Bremer%
}{%
{\protect \APACyear {2007}}%
}]{%
bremer2007varieties}
\APACinsertmetastar {%
bremer2007varieties}%
\begin{APACrefauthors}%
Bremer, M.%
\end{APACrefauthors}%
\unskip\
\newblock
\APACrefYearMonthDay{2007}{}{}.
\newblock
{\BBOQ}\APACrefatitle {Varieties of finitism} {Varieties of finitism}.{\BBCQ}
\newblock
\APACjournalVolNumPages{Metaphysica}{8}{2}{131--148}.
\PrintBackRefs{\CurrentBib}

\bibitem [\protect \citeauthoryear {%
Dummett%
}{%
Dummett%
}{%
{\protect \APACyear {1994}}%
}]{%
dummett1994mathematics}
\APACinsertmetastar {%
dummett1994mathematics}%
\begin{APACrefauthors}%
Dummett, M.%
\end{APACrefauthors}%
\unskip\
\newblock
\APACrefYearMonthDay{1994}{}{}.
\newblock
{\BBOQ}\APACrefatitle {What is Mathematics About?} {What is mathematics
  about?}{\BBCQ}
\newblock
\APACjournalVolNumPages{Mathematics and Mind}{}{}{11--26}.
\PrintBackRefs{\CurrentBib}

\bibitem [\protect \citeauthoryear {%
Enderton%
}{%
Enderton%
}{%
{\protect \APACyear {2001}}%
}]{%
enderton2001mathematical}
\APACinsertmetastar {%
enderton2001mathematical}%
\begin{APACrefauthors}%
Enderton, H\BPBI B.%
\end{APACrefauthors}%
\unskip\
\newblock
\APACrefYear{2001}.
\newblock
\APACrefbtitle {A mathematical introduction to logic} {A mathematical
  introduction to logic}.
\newblock
\APACaddressPublisher{}{Academic press}.
\PrintBackRefs{\CurrentBib}

\bibitem [\protect \citeauthoryear {%
Farmer%
}{%
Farmer%
}{%
{\protect \APACyear {2008}}%
}]{%
farmer2008seven}
\APACinsertmetastar {%
farmer2008seven}%
\begin{APACrefauthors}%
Farmer, W\BPBI M.%
\end{APACrefauthors}%
\unskip\
\newblock
\APACrefYearMonthDay{2008}{}{}.
\newblock
{\BBOQ}\APACrefatitle {The seven virtues of simple type theory} {The seven
  virtues of simple type theory}.{\BBCQ}
\newblock
\APACjournalVolNumPages{Journal of Applied Logic}{6}{3}{267--286}.
\PrintBackRefs{\CurrentBib}

\bibitem [\protect \citeauthoryear {%
Fletcher%
}{%
Fletcher%
}{%
{\protect \APACyear {1989}}%
}]{%
fletcher1989truth}
\APACinsertmetastar {%
fletcher1989truth}%
\begin{APACrefauthors}%
Fletcher, P.%
\end{APACrefauthors}%
\unskip\
\newblock
\APACrefYear{1989}.
\newblock
\APACrefbtitle {Truth, proof and infinity.} {Truth, proof and infinity.}
\newblock
\APACaddressPublisher{}{Kluver Academic Publishers, Netherlands}.
\PrintBackRefs{\CurrentBib}

\bibitem [\protect \citeauthoryear {%
Hallett%
}{%
Hallett%
}{%
{\protect \APACyear {1984}}%
}]{%
hallett1984cantorian}
\APACinsertmetastar {%
hallett1984cantorian}%
\begin{APACrefauthors}%
Hallett, M.%
\end{APACrefauthors}%
\unskip\
\newblock
\APACrefYear{1984}.
\newblock
\APACrefbtitle {Cantorian set theory and limitation of size} {Cantorian set
  theory and limitation of size}.
\newblock
\APACaddressPublisher{}{Clarendon Press Oxford}.
\PrintBackRefs{\CurrentBib}

\bibitem [\protect \citeauthoryear {%
Lavine%
}{%
Lavine%
}{%
{\protect \APACyear {2009}}%
}]{%
lavine2009understanding}
\APACinsertmetastar {%
lavine2009understanding}%
\begin{APACrefauthors}%
Lavine, S.%
\end{APACrefauthors}%
\unskip\
\newblock
\APACrefYear{2009}.
\newblock
\APACrefbtitle {Understanding the infinite} {Understanding the infinite}.
\newblock
\APACaddressPublisher{}{Harvard University Press}.
\PrintBackRefs{\CurrentBib}

\bibitem [\protect \citeauthoryear {%
Linnebo%
\ \BBA {} Shapiro%
}{%
Linnebo%
\ \BBA {} Shapiro%
}{%
{\protect \APACyear {2019}}%
}]{%
linnebo2019actual}
\APACinsertmetastar {%
linnebo2019actual}%
\begin{APACrefauthors}%
Linnebo, {\O}.%
\BCBT {}\ \BBA {} Shapiro, S.%
\end{APACrefauthors}%
\unskip\
\newblock
\APACrefYearMonthDay{2019}{}{}.
\newblock
{\BBOQ}\APACrefatitle {Actual and potential infinity} {Actual and potential
  infinity}.{\BBCQ}
\newblock
\APACjournalVolNumPages{No{\^u}s}{53}{1}{160--191}.
\PrintBackRefs{\CurrentBib}

\bibitem [\protect \citeauthoryear {%
Mostowski%
}{%
Mostowski%
}{%
{\protect \APACyear {2003}}%
}]{%
mostowski2003representing}
\APACinsertmetastar {%
mostowski2003representing}%
\begin{APACrefauthors}%
Mostowski, M.%
\end{APACrefauthors}%
\unskip\
\newblock
\APACrefYearMonthDay{2003}{}{}.
\newblock
{\BBOQ}\APACrefatitle {On representing semantics in finite models} {On
  representing semantics in finite models}.{\BBCQ}
\newblock
\BIn{} \APACrefbtitle {Philosophical dimensions of logic and science}
  {Philosophical dimensions of logic and science}\ (\BPGS\ 15--28).
\newblock
\APACaddressPublisher{}{Springer}.
\PrintBackRefs{\CurrentBib}

\bibitem [\protect \citeauthoryear {%
Mostowski%
\ \BBA {} Zdanowski%
}{%
Mostowski%
\ \BBA {} Zdanowski%
}{%
{\protect \APACyear {2005}}%
}]{%
mostowski2005fm}
\APACinsertmetastar {%
mostowski2005fm}%
\begin{APACrefauthors}%
Mostowski, M.%
\BCBT {}\ \BBA {} Zdanowski, K.%
\end{APACrefauthors}%
\unskip\
\newblock
\APACrefYearMonthDay{2005}{}{}.
\newblock
{\BBOQ}\APACrefatitle {FM-representability and beyond} {Fm-representability and
  beyond}.{\BBCQ}
\newblock
\BIn{} \APACrefbtitle {New Computational Paradigms} {New computational
  paradigms}\ (\BPGS\ 358--367).
\newblock
\APACaddressPublisher{}{Springer}.
\PrintBackRefs{\CurrentBib}

\bibitem [\protect \citeauthoryear {%
Mycielski%
}{%
Mycielski%
}{%
{\protect \APACyear {1986}}%
}]{%
Mycielski1986}
\APACinsertmetastar {%
Mycielski1986}%
\begin{APACrefauthors}%
Mycielski, J.%
\end{APACrefauthors}%
\unskip\
\newblock
\APACrefYearMonthDay{1986}{}{}.
\newblock
{\BBOQ}\APACrefatitle {Locally Finite Theories} {Locally finite
  theories}.{\BBCQ}
\newblock
\APACjournalVolNumPages{Journal of Symbolic Logic}{51}{1}{59--62}.
\PrintBackRefs{\CurrentBib}

\bibitem [\protect \citeauthoryear {%
Niebergall%
}{%
Niebergall%
}{%
{\protect \APACyear {2014}}%
}]{%
niebergall2014assumptions}
\APACinsertmetastar {%
niebergall2014assumptions}%
\begin{APACrefauthors}%
Niebergall, K\BHBI G.%
\end{APACrefauthors}%
\unskip\
\newblock
\APACrefYearMonthDay{2014}{}{}.
\newblock
{\BBOQ}\APACrefatitle {Assumptions of infinity} {Assumptions of
  infinity}.{\BBCQ}
\newblock
\APACjournalVolNumPages{Formalism and beyond: on the nature of mathematical
  discourse}{}{}{229--274}.
\PrintBackRefs{\CurrentBib}

\bibitem [\protect \citeauthoryear {%
Normann%
, Palmgren%
\BCBL {}\ \BBA {} Stoltenberg-Hansen%
}{%
Normann%
\ \protect \BOthers {.}}{%
{\protect \APACyear {1999}}%
}]{%
normann1999hyperfinite}
\APACinsertmetastar {%
normann1999hyperfinite}%
\begin{APACrefauthors}%
Normann, D.%
, Palmgren, E.%
\BCBL {}\ \BBA {} Stoltenberg-Hansen, V.%
\end{APACrefauthors}%
\unskip\
\newblock
\APACrefYearMonthDay{1999}{}{}.
\newblock
{\BBOQ}\APACrefatitle {Hyperfinite type structures} {Hyperfinite type
  structures}.{\BBCQ}
\newblock
\APACjournalVolNumPages{The Journal of Symbolic Logic}{64}{3}{1216--1242}.
\PrintBackRefs{\CurrentBib}

\bibitem [\protect \citeauthoryear {%
Parikh%
}{%
Parikh%
}{%
{\protect \APACyear {1971}}%
}]{%
parikh1971existence}
\APACinsertmetastar {%
parikh1971existence}%
\begin{APACrefauthors}%
Parikh, R.%
\end{APACrefauthors}%
\unskip\
\newblock
\APACrefYearMonthDay{1971}{}{}.
\newblock
{\BBOQ}\APACrefatitle {Existence and feasibility in arithmetic} {Existence and
  feasibility in arithmetic}.{\BBCQ}
\newblock
\APACjournalVolNumPages{The Journal of Symbolic Logic}{36}{03}{494--508}.
\PrintBackRefs{\CurrentBib}

\bibitem [\protect \citeauthoryear {%
Priest%
}{%
Priest%
}{%
{\protect \APACyear {2013}}%
}]{%
priest2013indefinite}
\APACinsertmetastar {%
priest2013indefinite}%
\begin{APACrefauthors}%
Priest, G.%
\end{APACrefauthors}%
\unskip\
\newblock
\APACrefYearMonthDay{2013}{}{}.
\newblock
{\BBOQ}\APACrefatitle {Indefinite extensibility—dialetheic style} {Indefinite
  extensibility—dialetheic style}.{\BBCQ}
\newblock
\APACjournalVolNumPages{Studia Logica}{101}{6}{1263--1275}.
\PrintBackRefs{\CurrentBib}

\bibitem [\protect \citeauthoryear {%
Sereno%
}{%
Sereno%
}{%
{\protect \APACyear {1998}}%
}]{%
sereno1998infinity}
\APACinsertmetastar {%
sereno1998infinity}%
\begin{APACrefauthors}%
Sereno, L\BPBI A.%
\end{APACrefauthors}%
\unskip\
\newblock
\APACrefYear{1998}.
\unskip\
\newblock
\APACrefbtitle {Infinity and experience} {Infinity and experience}\
  \APACtypeAddressSchool {\BUPhD}{}{}.
\unskip\
\newblock
\APACaddressSchool {}{Massachusetts Institute of Technology}.
\PrintBackRefs{\CurrentBib}

\bibitem [\protect \citeauthoryear {%
Shapiro%
\ \BBA {} Wright%
}{%
Shapiro%
\ \BBA {} Wright%
}{%
{\protect \APACyear {2006}}%
}]{%
shapiro2006all}
\APACinsertmetastar {%
shapiro2006all}%
\begin{APACrefauthors}%
Shapiro, S.%
\BCBT {}\ \BBA {} Wright, C.%
\end{APACrefauthors}%
\unskip\
\newblock
\APACrefYearMonthDay{2006}{}{}.
\newblock
{\BBOQ}\APACrefatitle {All things indefinitely extensible} {All things
  indefinitely extensible}.{\BBCQ}
\newblock
\APACjournalVolNumPages{Absolute generality}{}{}{255--304}.
\PrintBackRefs{\CurrentBib}

\bibitem [\protect \citeauthoryear {%
Tait%
}{%
Tait%
}{%
{\protect \APACyear {1981}}%
}]{%
tait1981finitism}
\APACinsertmetastar {%
tait1981finitism}%
\begin{APACrefauthors}%
Tait, W\BPBI W.%
\end{APACrefauthors}%
\unskip\
\newblock
\APACrefYearMonthDay{1981}{}{}.
\newblock
{\BBOQ}\APACrefatitle {Finitism} {Finitism}.{\BBCQ}
\newblock
\APACjournalVolNumPages{The Journal of Philosophy}{}{}{524--546}.
\PrintBackRefs{\CurrentBib}

\bibitem [\protect \citeauthoryear {%
Tait%
}{%
Tait%
}{%
{\protect \APACyear {2002}}%
}]{%
tait2002remarks}
\APACinsertmetastar {%
tait2002remarks}%
\begin{APACrefauthors}%
Tait, W\BPBI W.%
\end{APACrefauthors}%
\unskip\
\newblock
\APACrefYearMonthDay{2002}{}{}.
\newblock
{\BBOQ}\APACrefatitle {Remarks on finitism} {Remarks on finitism}.{\BBCQ}
\newblock
\APACjournalVolNumPages{Reflections on the Foundations of Mathematics. Essays
  in Honor of Solomon Feferman, LNL}{15}{}{}.
\PrintBackRefs{\CurrentBib}

\bibitem [\protect \citeauthoryear {%
Van~Bendegem%
}{%
Van~Bendegem%
}{%
{\protect \APACyear {1994}}%
}]{%
van1994strict}
\APACinsertmetastar {%
van1994strict}%
\begin{APACrefauthors}%
Van~Bendegem, J\BPBI P.%
\end{APACrefauthors}%
\unskip\
\newblock
\APACrefYearMonthDay{1994}{}{}.
\newblock
{\BBOQ}\APACrefatitle {Strict finitism as a viable alternative in the
  foundations of mathematics} {Strict finitism as a viable alternative in the
  foundations of mathematics}.{\BBCQ}
\newblock
\APACjournalVolNumPages{Logique et Analyse}{37}{145}{23--40}.
\PrintBackRefs{\CurrentBib}

\bibitem [\protect \citeauthoryear {%
Van~Bendegem%
}{%
Van~Bendegem%
}{%
{\protect \APACyear {1999}}%
}]{%
van1999largest}
\APACinsertmetastar {%
van1999largest}%
\begin{APACrefauthors}%
Van~Bendegem, J\BPBI P.%
\end{APACrefauthors}%
\unskip\
\newblock
\APACrefYearMonthDay{1999}{}{}.
\newblock
{\BBOQ}\APACrefatitle {Why the largest number imaginable is still a finite
  number} {Why the largest number imaginable is still a finite number}.{\BBCQ}
\newblock
\APACjournalVolNumPages{Logique et Analyse}{41}{}{161--162}.
\PrintBackRefs{\CurrentBib}

\end{thebibliography}

% The author information at the end of the article is printed automatically from the data supplied to \AuthorA, \AuthorB, and \AuthorC commands.

\end{document}